\documentclass[final,onefignum,onetabnum]{siamonline250211}
\usepackage{ulem}

\usepackage{braket,amsfonts,amssymb,multirow}

\usepackage{array}

\usepackage[caption=false]{subfig}
\usepackage{pgfplots}

\usepackage{graphicx,epstopdf,floatflt}

\usepackage{bold-extra}
\usepackage[most]{tcolorbox}

\colorlet{texcscolor}{blue!50!black}
\colorlet{texemcolor}{red!70!black}
\colorlet{texpreamble}{red!70!black}
\colorlet{codebackground}{black!25!white!25}

\lstdefinestyle{siamlatex}{%
  style=tcblatex,
  texcsstyle=*\color{texcscolor},
  texcsstyle=[2]\color{texemcolor},
  keywordstyle=[2]\color{texemcolor},
  moretexcs={cref,Cref,maketitle,mathcal,text,headers,email,url},
}

\tcbset{%
  colframe=black!75!white!75,
  coltitle=white,
  colback=codebackground, 
  colbacklower=white, 
  fonttitle=\bfseries,
  arc=0pt,outer arc=0pt,
  top=1pt,bottom=1pt,left=1mm,right=1mm,middle=1mm,boxsep=1mm,
  leftrule=0.3mm,rightrule=0.3mm,toprule=0.3mm,bottomrule=0.3mm,
  listing options={style=siamlatex}
}

\newtcblisting[use counter=example]{example}[2][]{%
  title={Example~\thetcbcounter: #2},#1}

\newtcbinputlisting[use counter=example]{\examplefile}[3][]{%
  title={Example~\thetcbcounter: #2},listing file={#3},#1}

\DeclareTotalTCBox{\code}{ v O{} }
{ 
  fontupper=\ttfamily\color{black},
  nobeforeafter,
  tcbox raise base,
  colback=codebackground,colframe=white,
  top=0pt,bottom=0pt,left=0mm,right=0mm,
  leftrule=0pt,rightrule=0pt,toprule=0mm,bottomrule=0mm,
  boxsep=0.5mm,
  #2}{#1}

\patchcmd\newpage{\vfil}{}{}{}
\begin{tcbverbatimwrite}{tmp_\jobname_header.tex}
\title{Bifurcation of Limit Cycles from a Fold-Fold Singularity in a Glacial Cycles Model
}

\author{Oleg Makarenkov\thanks{Department of Mathematical Sciences, University of Texas at Dallas, 75080 Richardson, USA, \email{makarenkov@utdallas.edu}.}
\and Esther Widiasih\thanks{Department of Mathematics, Natural and Health Sciences, University of Hawaii at West Oahu, 91-1001 Farrington Highway, Kapolei, HI 96707, USA, \email{widiasih@hawaii.edu}.}}

\end{tcbverbatimwrite}
\input{tmp_\jobname_header.tex}

\ifpdf
\hypersetup{ pdftitle={Guide to Using  SIAM'S \LaTeX\ Style} }
\fi

\newcommand{\eps}{\varepsilon}
\newcommand{\al}{\alpha}

\newsiamthm{cor}{Corollary}
\newsiamthm{df}{Definition}
\newsiamthm{lem}{Lemma}
\newsiamthm{rem}{Remark}
\newsiamthm{prop}{Proposition}

\usepackage{amsmath}

\begin{document}

\maketitle

\begin{tcbverbatimwrite}{tmp_\jobname_abstract.tex}
\begin{abstract}
We study the occurrence of limit cycles from a point on the discontinuity hyperplane $L$ between two smooth vector fields where the two vector fields both
point towards one another. Generically, such a point (called {\it switched equilibrium} in control) is asymptotically stable, but we consider the situation where the two vector fields become tangent to $L$ at the switched equilibrium under varying parameter making a degenerate fold-fold singularity. We prove that moving the parameter past such a singular value leads to the occurrence of an attracting limit cycle, which is exactly the dynamical mechanism we then discover in a conceptual model of glacial cycles. 
\end{abstract}


\end{tcbverbatimwrite}
\input{tmp_\jobname_abstract.tex}

\section{Background}

In the past 5 million years, planet Earth has gone through fascinating climate regimes. Started with the Pliocene Epoch about 5.3 to 2.8 million years before present (Myr BP) where the climate was 2-3$^o$C warmer than today, the planet cooled down during the later part of the epoch, and continued to do so in the following epoch, the Pleistocene, 2.5 Myr BP to about 12 thousand years ago. As the planet cooled down, Northern hemisphere ephemeral ice sheets turned into lasting glaciation during the late Pliocene, and eventually, in the Pleistocene Epoch, the glaciation period stayed longer and grew larger ice sheets, which  extended into the mid-latitude. At the last glacial maximum for example, the maximum ice extent went as far south as Wisconsin, which possibly was under a mile of ice \cite{stokes2017laurentide}. During the Pliocene Epoch, the period of glaciation follows the 41 kyr (1 kyr = 1000 years) period of obliquity signal of the Earth's orbit, having effects to the amount of incoming solar radiation (insolation). However, the glacial cycle transitioned from a 41 kyr periodicity to about 100 kyr at around the mid-Pleistocene, even though the insolation signal stays the same. The reason for such shift is still unknown, making the \emph{Mid Pleistocene Transition} (MPT) one of the most exciting questions in the scientific  community. Mathematically, this is an interesting problem as well, because understanding MPT requires an understanding of the non-linear feedbacks and responses of the system. \\

There is a hierarchy of models to describe the undeniably complex Earth's climate system. These models range from the most detailed general circulation models, to intermediate complexity models, to simpler conceptual climate models. Budyko's energy balance model (EBM) is a conceptual climate model of the planet's temperature distribution proposed in \cite{budyko1969effect}. The model was augmented with an equation of the ice line in \cite{widiasih2013dynamics} to capture the temperature-ice albedo feedback, which is one of the major climate feedbacks. The coupled temperature-ice line system was then coupled with an equation for the maximum ice extent, and the resulting \emph{glacial flip flop} model was shown to have a stable limit cycle \cite{wwhm2016periodic}. The system of three equations is capable of simulating the glacial cycles of the late Pleistocene, with or without the inclusion of time dependent orbital forcing. A current work shows that the glacial flip flop model admits the mid Pleistocene transition when a critical temperature parameter related to the cooling down of the planet is varied  \cite{widiasih2024mid}.  
 Indeed, a bifurcation in the glacial flip flop system relating to the cooling down of the planet gives way to the transition, and it is the goal of this paper to provide a mathematical proof and the analysis of its existence.

\vskip0.2cm

\noindent The glacial flip flop model can be broadly described as a Filippov system
\begin{equation}\label{modelbroad}
    \begin{array}{rcl}
    \dot x &=&f^i(x,y),\\
    \dot y &=&g^i(x,y,z,T^i),\\
    \dot z &=&h^i(x,y,z,T^i),
    \end{array}\quad i=\left\{\begin{array}{ll}
      +1,& z>0,\\
      -1,& z<0,
      \end{array}\right.
\end{equation}
where $T^-$ and $T^+$ are critical temperatures below which ice melts in glaciating and deglaciating states accordingly.
To catch the occurrence of Glacial cycles, we will view $T^-$ and $T^+$ as parameters and will obtain Glacial cycles as a bifurcation from a point as $T^-$ and $T^+$ vary.

\vskip0.2cm To outline the dynamical mechanism as for how attracting limit cycle in model (\ref{modelbroad}) bifurcates, two additional notions are needed. First, a point $(x,y,0)$ of the switching threshold
$$
  L=\left\{(x,y,z)\in\mathbb{R}^3:z=0\right\}
$$
is a switched equilibrium of (\ref{modelbroad}) if 
$$
  \left(\begin{array}{c}
    f^-\\
    g^-\\
    h^-
    \end{array}\right)(x,y,0,T^-)\uparrow\downarrow \left(\begin{array}{c}
    f^+\\
    g^+\\
    h^+
    \end{array}\right)(x,y,0,T^+),
$$
where symbol $\uparrow\downarrow$ says "oppositely directed".

In the next sections, we provide the details of the glacial flip flop model,  which is then followed by the proof of the bifurcation, and finally, a closing discussion.\\


\section{Model} 
\subsection{A Brief Recap of the Model Development and Analysis}
Mathematically speaking, the glacial flip flop model was built off three foundations: the construction of the integro-differential equation of the temperature distribution known as the Budyko EBM \cite{budyko1969effect}, the inclusion of the ice/ albedo line equation \cite{widiasih2013dynamics}, and the dimension reduction in \cite{mcgwid2014simplification}. The infinite dimensional system of temperature distribution $T(y)$ and ice/ albedo line latitudinal position $\eta$ in \cite{mcgwid2014simplification} is reduced to a system of 2 dimensional ordinary differential equations of $w$ and $\eta$, where $w$ is a translate of a global mean temperature. 

The $w -\eta$ system is then coupled with an equation of the ice extent latitudinal position $\xi$ in \cite{wwhm2016periodic} as non-smooth system of differential equations. The analysis in \cite{wwhm2016periodic} included the details of the modeling approach as well as the existence proof of the limit cycle,  applying the contraction mapping principle. The work in this paper is to provide an alternative proof that the bifurcation comes from a fold-fold singularity of the critical temperature parameter.


Bifurcation of attracting cycles from a point of a discontinuity manifold $L$ is addressed in  \cite{Simpson,Kuepper2,Kuepper1} when the point is an equilibrium for each of the vector fields on the two sides of $L.$ Bifurcation of limit cycles from a  fold-fold singularity of $L$ was considered in  \cite{MakSIAMfold-fold,MakFranklin} in two-dimensional settings. The fold-fold singularity in 3D case is also called {\it Teixeira singularity} \cite{Teixeira} and it is known to be structurally stable. In other words, fold-fold singularity of $L$ in 3D can produce limit cycles only when the singularity is degenerate, i.e. when the vector fields on the two sides of $L$ are not just tangent one another, but pointing opposite one another (known as {\it switched equilibrium} in control \cite{Bolzern}). Bifurcation of limit cycles from a degenerate fold-fold singularity for particular quadratic vector fields is investigated in \cite{Cristiano}.  

\subsection{Nonsmooth System of the Glacial Cycles}

In this work, the glacial cycle is modeled through 3 climate variables: the temperature, $w$, the albedo/ snow line $\eta$, and the ice sheet extent, $\xi$ \cite{wwhm2016periodic}. The \emph{temperature-snow line-ice sheet extent} model can be written as the following non-smooth dynamical systems:

\begin{align}
\dot{w} &=\tau \left( w-F(\eta) \right) \notag\\
\dot{\eta} &= \rho \left( w- G^{\pm}(\eta) \right) \hspace{1cm} \text{ when } (w, \eta,\xi) \in V^{\pm} \label{eqn_wetaxidot} \\
\dot{\xi} &= \epsilon H^{\pm} \left(\eta, \xi \right) \hspace{1.8cm} \text{ when } (w, \eta,\xi) \in V^{\pm} \notag
\end{align}

\begin{align}
&F(\eta)=\frac{1}{B}\left(Q(1-\al_0)-A+C \frac{Q}{B+C}(\al_2-\al_1)(\eta-\textstyle{\frac{1}{2}}+s_2P_2(\eta))\right).\\
&G^{\pm}(\eta) = L s_2(1-\alpha_0)p_2(\eta)+T^{\pm} \\
&H^{\pm}(\eta, \xi) = b^{\pm}(\eta-\xi)-a(1-\eta)
\end{align}
 where each equation is defined in its respective region: 
\begin{align} 
&V^{-} := \{ \, (w,\eta,\xi) : b(\eta-\xi)-a(1-\eta)<0 \, \} \text{ glaciating state }\\
&V^{+} := \{ \, (w,\eta,\xi) : b(\eta-\xi)-a(1-\eta)>0 \, \} \text{ deglaciating state}
\end{align}
with the 2nd order Legendre function $p(\eta)$ and its indefinite integral $P(\eta)$  being:\\
\begin{align*}
p(\eta) &= \frac{1}{2} (3 \eta^2-1)\\
P(\eta) &= \frac{1}{2}(\eta^3 - \eta)\\
\end{align*}

The parameter values as used in \cite{wwhm2016periodic} are reproduced in Table \ref{table1}

\begin{center}{\small}
\begin{tabular}{|ccc|ccc|}\hline
Parameter & Value & Units & Parameter & Value & Units\\ \hline
$Q = \text{constant}$  & 343 (autonomous) & Wm$^{-2}$ & $T_c^-$ & $ -5.5 $ & $^\circ$C\\ 
$Q=Q(t) $  & as in Eq. \eqref{eqn_Q} & Wm$^{-2}$ & $T_c^+$ & -10  & $^\circ$C\\ 
$s_2 = \text{constant}$  & -0.482 (autonomous) & dim'less & $s_2(t)$ & as in Eq. \eqref{eqn_s2} & dim'less\\ 
$A$ & 202 & Wm$^{-2}$ & $b^{+}$ & 1.5 & dim'less\\ 
$B$ & 1.9 & Wm$^{-2}(^\circ\mbox{C})^{-1}$ & $b$ & 1.75 & dim'less\\
$C$ & 3.04 & Wm$^{-2}(^\circ\mbox{C})^{-1}$ & $b^{-}$ & 5 & dim'less\\
$\al_1$ & 0.32 & dim'less & $a$ & 1.45 & dim'less\\
$\al_2$ & 0.62 & dim'less & $\tau$ & 7 & seconds$^{-1}$\\ 
$T_1$ & $-10$ & $^\circ$C & $\rho$ & 0.05 & seconds$^{-1} (^\circ\mbox{C})^{-1}$\\
\quad & \quad & \quad & $\epsilon$ & 0.03 & \text{seconds}$^{-1}$\\
 \hline
\end{tabular}\label{table1}
\end{center}
\vspace{0.1in}
{\bf Table \ref{table1}. Values of the parameters used to simulate system \eqref{eqn_wetaxidot} as reproduced from \cite{wwhm2016periodic}}\\

The Milankovitch orbital forcing enters through $Q$ and $s_2$ as computed in \cite{mcglehman2012} and \cite{nadeau2017simple}, as time dependent functions of the Earth's \emph{eccentricity}, $e = e(t)$ and \emph{obliquity}, $\beta = \beta(t)$, using the computational result in \cite{laskar2004long}:

\begin{align} 
&Q(t) = Q(e) =  \frac{Q_0}{\sqrt{1-e^2}} \quad \text{ with } Q_0 = 343 \text{ Wm}^{-2} \label{eqn_Q}\\
&s_2(t) = s_2( \beta) = \frac{5}{16} \left( 3 \cos^2(\beta) - 1 \right) \label{eqn_s2} 
\end{align}

\section{Analysis}

To prove the existence of an attractive limit cycle in the flip flop model, we rewrite the system for convenience:
\begin{equation}\label{sys0}
   \begin{array}{rcl}
       \dot w & = & -\tau(w-F(\eta)),\\
       \dot \eta & = & \rho (w+L s_2 (1-\alpha_0)p_2(\eta))-\rho (T_c^\pm+\eps \bar T^\pm_c),\\
       \dot \xi &=&\epsilon(b^\pm(\eta-\xi)-a(1-\eta)),
    \end{array}
\end{equation}
we introduce a small parameter $\eps>0$ as
$$
\begin{array}{lllllll}
   b^-=b_0+\eps\bar b_0, & & b^+=b_1+\eps \bar b_1, && \ b=1.75, & b_0=1.5, & b_1=5,\\
    T_c^-=T^-+\eps \overline T^-, & & T_c^+=T^++\eps \overline T^+,  &&   T_c^- = -10, \  T_c^+ \approx T_c^-
\end{array}
$$
The work in \cite{wwhm2016periodic} suggests that minus or plus in (\ref{sys0}) take place according to whether $(b+a)\eta-a-\xi b$ is negative or positive. We therefore, need to study the occurence of limit cycles in the system
\begin{equation}\label{sys1}
   \begin{array}{rcl}
       \dot x & = & f(x,y),\\
       \dot y & = & g(x,y)+g^i(\eps),\\
       \dot z & = & h^i(y,z,\eps),
    \end{array}\quad 
    q=\left\{\begin{array}{ll}
                -1, &{\rm if} \ H(y,z)<0,\\
                +1, & {\rm if} \ H(y,z)>0,
                \end{array}\right.
    \quad H(y,z)=(b+a) y-a-bz.
\end{equation}
$$
   l_0\in \{(x,y,z): H(y,z)=0\}=:L.
$$
We plan to prove bifurcation of a stable limit cycle from a point $l_0.$ 

\vskip0.2cm

\noindent In what follows, we denote by $t\mapsto\left(X^i(t,\bar x,\bar y,\bar z,\eps),Y^i(t,\bar x,\bar y,\bar z,\eps),Z^i(t,\bar x,\bar y,\bar z,\eps)\right)^T$ the general solution of the $i$-th subsystem of (\ref{sys1}) with the initial condition $x(0)=\bar x$, $y(0)=\bar y,$ $z(0)=\bar z.$ 
Throughout the paper we allow ourselves to identify vector lines and vector columns where it doesn't lead to a confusion.

\vskip0.2cm

\noindent The following notations are required to formulate the theorem. They will also be used throughout the proof.
 \begin{eqnarray*}
     \bar\alpha^i &=&-2\dfrac{(a+b)g_x(x_0,y_0)}{D^i},\\
\bar\gamma^i&=&-2\dfrac{(a+b)(g^i)'(0)-bh^i_\eps(y_0,z_0,0)}{D^i},\\      
D^i&=&bg_y(x_0,y_0)-\dfrac{b^2}{a+b}h^i_y(y_0,z_0,0)-bh_z^i(y_0,z_0,0)
\end{eqnarray*}
 $$  \zeta(z)=\frac{a+bz}{a+b},\quad \zeta'(0)=\dfrac{b}{a+b}$$  
  \begin{eqnarray*}
     \bar {\bar k}^i &=& \left.\dfrac{a+b}{6}\left(g'_x(x_0,y_0)f'_y(x_0,y_0)+g{}'_y{}'_y(x_0,y_0)(g(x_0,y_0)+ g^i(0))+(g'_y(x_0,y_0))^2\right)\right.-\\
     && {\color{black}-\dfrac{b}{6}\left(h^i{}'_y{}'_y(y_0,z_0,0)+\dfrac{a+b}{b}h^i{}'_y{}'_z(y_0,z_0,0)+\dfrac{a+b}{b}h^i{}'_z{}'_y(y_0,z_0,0)+\dfrac{(a+b)^2}{b^2}h^i{}'_z{}'_z(y_0,z_0,0)\right)-}\\
     &&\left.-\dfrac{b}{6}\left(h^i{}'_y(y_0,z_0,0)g'_y(x_0,y_0)+h^i{}'_z(y_0,z_0,0)h^i{}'_y(y_0,z_0,0)+\dfrac{a+b}{b}h^i{}'_z(y_0,z_0,0)^2\right)\right.
\end{eqnarray*}
\begin{eqnarray*}
    \tilde k^i &=&\dfrac{a+b}{2}\left(g'_x(x_0,y_0)f'_y(x_0,y_0)+g{}'_y{}'_y(x_0,y_0)(g(x_0,y_0)+\tilde g^i(0))+g'_y(x_0,y_0)g'_y(x_0,y_0)\right)\zeta'(0)-\\
     && {\color{black}-\dfrac{b}{2}\left(h^i{}'_y{}'_y(y_0,z_0,0)+\dfrac{a+b}{b}h^i{}'_z{}'_y(y_0,z_0,0)\right)(g(x_0,y_0)+ g^i(0))\zeta'(0)-}\\
     &&-\dfrac{b}{2}\left(h^i{}'_y(y_0,z_0,0)g'_y(x_0,y_0)+h^i{}'_z(y_0,z_0,0)h^i{}'_y(y_0,z_0,0)\right)\zeta'(0)-\\
     &&{\color{black}-\dfrac{b}{2}\left(h^i{}'_y{}'_z(y_0,z_0,0)+\dfrac{a+b}{b}h^i{}'_z{}'_z(y_0,z_0,0)\right)(g(x_0,y_0)+ g^i(0))}{-\dfrac{b}{2}(h^i{}'_z(y_0,z_0,0))^2},
     \end{eqnarray*}
\begin{eqnarray*}
     \hat k^i&=&\dfrac{a+b}{2}g{}'_y{}'_y(x_0,y_0)\zeta'(0)^2{\color{black}-b\left(\dfrac{1}{2}h^i{}'_y{}'_y(y_0,z_0,0)\zeta'(0)^2+h^i{}'_y{}'_z(y_0,z_0,0)\zeta^i(0)+\dfrac{1}{2}h^i{}'_z{}'_z(y_0,z_0,0)\right)}\\ 
     \bar\eta^i&=&-4\dfrac{1}{D^i}\left(4\dfrac{b}{a+b}\bar {\bar k}^i \dfrac{1}{h^i(y_0,z_0,0)}-\tilde k^i\dfrac{2}{h^i(y_0,z_0,0)}+\hat k^i\right),
\end{eqnarray*}
\begin{eqnarray*}
     B^i&=&\dfrac{1}{2}f'_y(x_0,y_0)\dfrac{b}{a+b}\bar\eta^i+4\dfrac{b}{a+b}\left(\dfrac{1}{6}f'_x(x_0,y_0)g'_y(x_0,y_0)+\right.\\
     && \left.+\dfrac{1}{6}f{}'_y{}'_y(x_0,y_0)(g(x_0,y_0)+g^i(0))+\dfrac{1}{6}f'_y(x_0,y_0)g'_y(x_0,y_0)\right)\dfrac{1}{h^i(y_0,z_0,0)}+\\
     &&-\left(f_x(x_0,y_0)f_y(x_0,y_0)+f_y(x_0,y_0)g_y(x_0,y_0)\right)\dfrac{b}{a+b}\cdot\dfrac{1}{h^i(y_0,z_0,0)}-\dfrac{1}{2}f_{yy}(x_0,y_0)\dfrac{b^2}{(a+b)^2},
\end{eqnarray*}
\begin{eqnarray*}
       k&=&-\dfrac{\bar \eta^--\bar\eta^+ }{\bar\alpha^--\bar\alpha^+},\\   
  m&=&-\dfrac{\bar\gamma^--\bar\gamma^+}{\bar\alpha^--\bar\alpha^+}, \\
A^i&=&f_x(x_0,y_0)+
\dfrac{1}{2}f_y(x_0,y_0)(g(x_0,y_0)+g^i(0)), \\   
 K&=&2\dfrac{kA^++B^+}{h^+(y_0,z_0,0)}-2\dfrac{kA^-+B^-}{h^-(y_0,z_0,0)}\\
     C^i&=& \dfrac{1}{2}f_y(x_0,y_0)\dfrac{b}{a+b}\bar\gamma^i.\\
M&=&(mA^-+C^-)\beta^--(mA^++C^+)\beta^+.
  \end{eqnarray*}
Note, $\nabla H(y,z)$ doesn't depend on $(y,z)$. Specifically,
 $$
 \nabla H(y,z)=(a+b,-b)=:\nabla H
 $$

\begin{theorem} Let 
\begin{equation}\label{onemore}
(x_0,y_0,z_0)\in L
\end{equation}  be such a point that
\begin{equation}\label{4}
   f(x_0,y_0)=0 \quad \mbox{\rm (border collision in {\it x}  variable)}
\end{equation}
and
\begin{equation}\label{as1}
\nabla H\left(\begin{array}{c} g(x_0,y_0)+ g^i(0)\\h^i(y_0,z_0,0)\end{array}\right)=0,\quad i\in\{-1,+1\} \quad\mbox{\rm (fold-fold in ({\it y,z})  variables)}.
\end{equation}
Assume that 
\begin{equation}\label{as3}
\bar\alpha^- -\bar\alpha^+\not=0, 
\end{equation}
\begin{equation}\label{as4}
KM<0\quad \mbox{\rm (bifurcation of 2 fixed points occurs)}.
\end{equation}
Finally, assume that 
\begin{eqnarray}
 &&  h^+(y_0,z_0,0)\cdot(\bar\eta^+-\bar\eta^-)<0,\qquad\hskip0.04cm \mbox{\rm (stability of 1/2 fixed points)} \label{stab}  \\
&&   h^+(y_0,z_0,0)h^-(y_0,z_0,0)<0 \label{h+h-},\qquad \mbox{\rm (positivity of time map of 1/2 fixed points)}\label{timemap}\\
&&bh^+(y_0,z_0,0)\cdot \left(g_y(x_0,y_0)-\dfrac{b}{a+b}h^+_y(y_0,z_0,0)-h^+_z(y_0,z_0,0)\right)<0,\label{directions1}\\
&&bh^+(y_0,z_0,0)\cdot \left(g_y(x_0,y_0)-\dfrac{b}{a+b}h^-_y(y_0,z_0,0)-h^-_z(y_0,z_0,0)\right)<0.\label{directions2}\\
&& \mbox{\rm (conditions for the limit cycle to be real rather than virtual)}
\end{eqnarray}
Then, for all $\eps>0$ sufficiently small, system (\ref{sys1}) admits a unique attractive limit cycle in a small neighborhood of  $(x_0,y_0,z_0)$ that shrinks to $(x_0,y_0,z_0)$ as $\eps\to 0.$ The period of the cycle equals
\begin{equation}\label{Tformula}
   T=\dfrac{2}{|h^-(y_0,z_0,0)|}\sqrt{-\eps M/K}+\dfrac{2}{|h^+(y_0,z_0,0)|}\sqrt{-\eps M/K}+
   O(\eps).
\end{equation}
\end{theorem}

\noindent {\bf Proof.} 
 {\bf Step 1: }{\it Expanding the time map $T^i_\eps(x,z)$.}  We will find $T^i_\eps(x,z)$ as a solution of the equation
\begin{equation}\label{temp1}
\begin{array}{l}
   H(Y^i(T,x,\zeta(z),z,\eps),Z^i(T,x,\zeta(z),z,\eps))=\\
   \ \ =(a+b)Y^i(T,x,\zeta(z),z,\eps)-a-bZ^i(T,x,\zeta(z),z,\eps)=0.
   \end{array}
\end{equation}
To do this, we expand $Y^i(T,x,\zeta(z),z,\eps)$ and $Z^i(T,x,\zeta(z),z,\eps)$ in Taylor's series near $(T,x,z,\eps)=(0,x_0,z_0,0)$. We write down an expansion for $Z^i$ only as $Y^i$ expands analogously.
$$
    \renewcommand{\arraystretch}{1.5}
\begin{array}{l}    Z(T,x,\zeta(z),z,\eps)=Z(0,x_0,\zeta(z_0),z_0,0)+Z^i{}'_t(l_{00})T+Z^i{}'_y(l_{00})\zeta'(0)(z-z_0)+Z^i{}'_z(l_{00})(z-z_0)
+\\
\qquad+\left[{\color{black}\dfrac{1}{2}}Z^i{}'_t{}'_t(l_{00})T+Z^i{}'_t{}'_x(l_{00})(x-x_0)+\left(Z^i{}'_t{}'_y(l_{00})\zeta'(0)+Z^i{}'_t{}'_z(l_{00})\right)(z-z_0)+Z^i{}'_t{}'_\eps(l_{00})\eps\right]T+\\
\qquad +{\color{black} \left[\dfrac{1}{6}Z^i{}'_t{}'_t{}'_t(l_{00})T^2+\dfrac{1}{2}\left(Z^i{}'_t{}'_t{}'_y(l_{00})\zeta'(0)+Z^i{}'_t{}'_t{}'_z(l_{00})\right)T(z-z_0)+\right.} \\ 
\qquad\qquad{\color{black}\left.+\left(\dfrac{1}{2}Z^i{}'_t{}'_y{}'_y(l_{00})\left(\zeta'(0)\right)^2+Z^i{}'_t{}'_y{}'_z(l_{00})\zeta'(0)+\dfrac{1}{2}Z^i{}'_t{}'_z{}'_z(l_{00})\right)(z-z_0)^2\right]T}+\mbox{remaining terms.}
\end{array}
$$
Using these expansions along with assumptions (\ref{4}) and (\ref{as1}) we can  rewrite (\ref{temp1}) as 
$$
    \renewcommand{\arraystretch}{1.5}
\begin{array}{l}  
   \bar m^i T+\tilde m^i(x-x_0)+\hat m^i(z-z_0)+m^i\eps+\\
\qquad\qquad+\bar k^i T^2+\tilde k^i T(z-z_0)+\hat k^i(z-z_0)^2+
    \mbox{remaining terms}=0,
\end{array}
$$
where
\begin{eqnarray*}
   \bar m^i&=& \dfrac{1}{2}(a+b)Y^i{}'_t{}'_t(l_{00})-\dfrac{1}{2}bZ^i{}'_t{}'_t(l_{00}),\\
   \tilde m^i&=& (a+b)Y^i{}'_t{}'_x(l_{00})-bZ^i{}'_t{}'_x(l_{00}),\\
   \hat m^i &=&(a+b)(Y^i{}'_t{}'_y(l_{00})\zeta'(0)+Y^i{}'_t{}'_z(l_{00}))-b(Z^i{}'_t{}'_y(l_{00})\zeta'(0)+Z^i{}'_t{}'_z(l_{00})),\\
   m^i&=&(a+b)Y^i{}'_t{}'_\eps(l_{00})-bZ^i{}'_t{}'_\eps(l_{00}),\\
   \bar k^i & = & \dfrac{1}{6}(a+b)Y^i{}'_t{}'_t{}'_t(l_{00})-\dfrac{1}{6}bZ^i{}'_t{}'_t{}'_t(l_{00}),\\
   \tilde k^i &=& \dfrac{1}{2}(a+b)(Y^i{}'_t{}'_t{}'_y(l_{00})\zeta'(0)+Y^i{}'_t{}'_t{}'_z(l_{00}))-\dfrac{1}{2}b(Z^i{}'_t{}'_t{}'_y(l_{00})\zeta'(0)+Z^i{}'_t{}'_t{}'_z(l_{00})),\\
   \hat k^i &=& (a+b)\left(\dfrac{1}{2}Y^i{}'_t{}'_y{}'_y(l_{00})(\zeta'(0))^2+Y^i{}'_t{}'_y{}'_z(l_{00})+\dfrac{1}{2}Y^i{}'_t{}'_z{}'_z(l_{00})\right)-\\
   && - b\left(\dfrac{1}{2}Z^i{}'_t{}'_y{}'_y(l_{00})(\zeta'(0))^2+Z^i{}'_t{}'_y{}'_z(l_{00})+\dfrac{1}{2}Z^i{}'_t{}'_z{}'_z(l_{00})\right).
\end{eqnarray*}
Therefore, we can compute $T^i_\eps$ as
\begin{equation}\label{Ti}
   T^i_\eps(x,z)= \alpha^i(x-x_0)+\beta^i(z-z_0)+\gamma^i\eps+\eta^i(z-z_0)^2+\mbox{remaining terms},
\end{equation}
where
$$
   \alpha^i=-\dfrac{\tilde m^i}{\bar m^i},\ \ \beta^i=-\dfrac{\hat m^i}{\bar m^i},\ \ \gamma^i=-\dfrac{m^i}{\bar m^i}, \ \ \eta^i=-\dfrac{\bar k^i(\beta^i)^2+\tilde k^i\beta^i+\hat k^i}{\bar m^i}.
$$

\noindent {\bf Step 2:} {\it The Poincare map $P_\eps$ and its fixed points.}  Expanding $X^i(T,x,\zeta(z),z,\eps)$ about $(T,x,z,\eps)=(0,x_0,z_0,0)$ we get
$$
 \renewcommand{\arraystretch}{2}
\begin{array}{l}
X^i(T,x,\zeta(z),z,\eps)=x_0+X^i{}'_t(l_{00})T+X^i{}'_x(l_{00})(x-x_0)+\\
+\left[{\color{black}\dfrac{1}{2}}X^i{}'_t{}'_t(l_{00})T+X^i{}'_t{}'_x(l_{00})(x-x_0)+X^i{}'_t{}'_y(l_{00})\zeta'(0)(z-z_0)+X^i{}'_t{}'_z(l_{00})(z-z_0)+X^i{}'_t{}'_\eps(l_{00})\eps\right]T+\\
+\left[\dfrac{1}{6}X{}'_t{}'_t{}'_t(l_{00})T^2+\dfrac{1}{2}X^i{}'_t{}'_t{}'_y(l_{00})T\zeta'(0)(z-z_0)+\dfrac{1}{2}X^i{}'_t{}'_t{}'_z(l_{00})T(z-z_0)+\right.\\+\left.\left(\dfrac{1}{2}X{}'_t{}'_y{}'_y(l_{00})\zeta'(0)^2+X{}'_t{}'_y{}'_z(l_{00})\zeta'(0)+\dfrac{1}{2}X{}'_t{}'_z{}'_z(l_{00})\right)(z-z_0)^2\right]T
+\mbox{remaining terms},
\end{array}
$$
and similarly for $Z^i(T,x,\zeta(z),z,\eps).$ 
Therefore, for the map
$$
  P^i_\eps(x,z)=\left(\begin{array}{c} X^i(T^i_\eps(x,z),x,\zeta(z),z,\eps) \\ Z^i(T^i_\eps(x,z),x,\zeta(z),z,\eps)\end{array}\right)
$$
we have 
$$
\renewcommand{\arraystretch}{2}
\begin{array}{l}
P^i_\eps(x,z)=\left(\begin{array}{l} x+\left[{\color{black}\dfrac{1}{2}}X^i{}'_t{}'_t(l_{00})T^i(x,z,\eps)+X^i{}'_t{}'_x(l_{00})(x-x_0)+X^i{}'_t{}'_y(l_{00})\zeta'(0)(z-z_0)+\right. \\ \qquad +\left.\dfrac{1}{6}X^i{}'_t{}'_t{}'_t(l_{00})T^i(x,z,\eps)^2+\dfrac{1}{2}X^i{}'_t{}'_t{}'_y(l_{00})T^i(x,z,\eps)\zeta'(0)(z-z_0)+\right. \\ \qquad +\left.\dfrac{1}{2}X^i{}'_t{}'_y{}'_y(l_{00})\zeta'(0)^2(z-z_0)^2\right] T^i(x,z,\eps) \\ z+Z^i{}'_t(l_{00})T^i(x,z,\eps)+\left[{\color{black}\dfrac{1}{2}}
    Z^i{}'_t{}'_t(l_{00})T^i_\eps(x,z)+Z^i{}'_t{}'_y(l_{00})\zeta'(0)(z-z_0)+\right.\\
    \quad\left.\hskip4cm+Z^i{}'_t{}'_z(l_{00})(z-z_0)\right]T^i(x,y,\eps)\end{array}\right)+\\
 \qquad\qquad\quad +\  \mbox{remaining terms}.
\end{array}
$$
Now we cancel some terms by observing that (\ref{as1}) yields
\begin{equation}\label{yields}
 \beta^i=-\dfrac{2}{h^i(y_0,z_0,0)},
 \end{equation}
 from where
 \begin{eqnarray*}
 \dfrac{1}{2}X^i{}'_t{}'_t(l_{00})\beta^i+X^i{}'_t{}'_y(l_{00})\zeta'(0)&=&0,\\
 \dfrac{1}{2}Z^i{}'_t{}'_t(l_{00})\beta^i+Z^i{}'_t{}'_y(l_{00})\zeta'(0)+Z^i{}'_t{}'_z(l_{00})&=&0,
\end{eqnarray*}
so that 
  $P_\eps^i$ can be rewritten in the form
\begin{eqnarray}
  P^i_\eps(x,z)&=&\left(\begin{array}{c} x\\ z\end{array}\right)+\left(\begin{array}{c}
  A^i(x-x_0)+B^i (z-z_0)^2+C^i\eps \\ h^i(y_0,z_0,0)\end{array}\right)\cdot \nonumber \\
  &&
  \cdot \left((\alpha^i,\beta^i)\left(\begin{array}{c}
  x-x_0\\ z-z_0\end{array}\right)+\eta^i(z-z_0)^2+\gamma^i\eps\right) +\,\mbox{remaining terms},\label{P}
\end{eqnarray}
where
\begin{eqnarray*}
   A^i&=&\dfrac{1}{2}X^i_{tt}(l_{00})+X^i{}'_t{}'_x(l_{00}),\\ 
B^i&=&\dfrac{1}{2}X^i{}'_t{}'_t(l_{00})\eta^i+\dfrac{1}{6}X^i{}'_t{}'_t{}'_t(l_{00})(\beta^i)^2+\dfrac{1}{2}X^i{}'_t{}'_t{}'_y(l_{00})\zeta'(0)\beta^i+\dfrac{1}{2}X^i{}'_t{}'_y{}'_y(l_{00})\zeta'(0)^2,\\ C^i&=&\dfrac{1}{2}X^i{}'_t{}'_t(l_{00})\gamma^i.
\end{eqnarray*}
Since the way $P_\eps^i$ is introduced implies that $P_\eps^i=\left(P_\eps^i\right)^{-1}$ (i.e. $P_\eps^i$ is an involution) we can find fixed points of the map
$$
   P_\eps(u)=P_{\eps}^-(P_\eps^+(u))
$$
by solving the equation
\begin{equation}\label{PP}
   P_{\eps}^-(u)=P_\eps^+(u).
\end{equation}
We will first solve the second equation of (\ref{PP}) and find $x(z,\eps)$. This will be found uniquely. Then we will plug the result into the first equation of (\ref{PP}) and find $z(\eps).$ The later will have two solutions which correspond to the two points where the cycle intersects the cross-section $L$.

\vskip0.2cm

\noindent Letting
$$
   \Phi(x,z,\eps)=\left[P_\eps^-(x,z)-P_\eps^+(x,z)\right]_2,
$$
we compute 
\begin{eqnarray*}
   \Phi'_{x}(x_0,z_0,0)&=&\alpha^- h^-(y_0,z_0,0)-\alpha^+h^+(y_0,z_0,0),\\
\Phi'_{z}(x_0,z_0,0)&=&\beta^- h^-(y_0,z_0,0)-\beta^+h^+(y_0,z_0,0)=0,   \\
\Phi'_z{}'_z(x_0,z_0,0)&=&2\eta^-h^-(y_0,z_0,0)-2\eta^+ h^+(y_0,z_0,0),\\
\Phi'_{\eps}(x_0,z_0,0)&=&\gamma^-h^-(y_0,z_0,0)-\gamma^+h^+(y_0,z_0,0),
\end{eqnarray*}
Using assumption (\ref{as3}) we  apply the Implicit Function Theorem and solve $\Phi(x,z,\eps)=0$ in $x$ when $(z,\eps)$ is near $(z_0,0).$ The Implicit Function Theorem gives
$$
   (x'_z,x'_\eps)(z_0,0)=-\dfrac{1}{\Phi'_{x}(x_0,z_0,0)}\left(0,\Phi'_{\eps}(x_0,z_0,0)\right),\quad x'_z{}'_z(z_0,0)=-\dfrac{1}{\Phi'_{x}(x_0,z_0,0)}\Phi'_z{}'_z(x_0,z_0,0).
$$
Therefore,
\begin{equation}\label{u1}
   x(z,\eps)=x_0+k (z-z_0)^2+
   m\eps+\mbox{remaining terms}.
\end{equation}
Plugging expression (\ref{u1}) into the first line of (\ref{PP}) we obtain the following equation for $z-z_0$
$$
   J(z-z_0)^4+K (z-z_0)^3+L\eps(z-z_0)^2+M\eps(z-z_0)+N\eps^2+\mbox{remaining terms}=0.
$$
The change of the variables
$$
   z-z_0=\eps^{1/2}p
$$
yields
$$
   J\eps^{2}p^4+K\eps^{3/2}p^3+L\eps^2 p^2+M\eps^{3/2}p+N\eps^{2}+\mbox{remaining terms}=0,
$$   
or
$$
   J\eps^{1/2}p^4+K p^3+L\eps^{1/2} p^2+M p+N\eps^{1/2}+\mbox{remaining terms}=0,
$$   
where we will need formulas for only some of $J,$ $K,$ $L,$ $M,$ and $N$, that we give a bit later.
Assumption (\ref{as4}) ensures that this quadratic polynomial admits two solutions $\underline z-z_0$ and $\overline z-z_0$. Combining formulas for solutions of the quadratic polynomial with with formula (\ref{u1}), we conclude that the Poincare map $P_\eps$ possesses the following two fixed points around $(x_0,z_0)$ for $\eps>0$ sufficiently small
\begin{equation}\label{xz}
\begin{array}{l}
   \begin{array}{r}\left(\begin{array}{c}\underline{x}(\eps) \\ \underline{z}(\eps)\end{array}\right)=\left(\begin{array}{c}
   x_0\\ z_0\end{array}\right)+\left(\begin{array}{c}\eps(-kM/K+m) \\ -\eps^{1/2}\sqrt{-M/K}\end{array}\right)+\mbox{remaining terms},
\end{array}\\ \begin{array}{r}\left(\begin{array}{c}\overline{x}(\eps) \\ \overline{z}(\eps)\end{array}\right)=\left(\begin{array}{c}
   x_0\\ z_0\end{array}\right)+\left(\begin{array}{c}\eps(-kM/K+m) \\ \eps^{1/2}\sqrt{-M/K}\end{array}\right)+\mbox{remaining terms},\end{array}
\end{array}
\end{equation}
where
\begin{eqnarray*}
K&=&(kA^-+B^-)\beta^--(kA^++B^+)\beta^+,\\
M&=&(mA^-+C^-)\beta^--(mA^++C^+)\beta^+.
\end{eqnarray*}

\noindent {\bf Step 3:} {\it Stability of fixed points $(\underline x(\eps),\underline z(\eps))$ and $(\overline x(\eps),\overline z(\eps))$.} The fixed point $(\underline x(\eps),\underline z(\eps))$ is stable, if the eigenvalues of the matrix 
\begin{equation}\label{compo}
   (P_\eps)'(\underline x(\eps),\underline z(\eps))=(P_\eps^-)'(\underline{\underline{x}}(\eps),\underline{\underline{z}}(\eps))\circ(P_\eps^+)'(\underline x(\eps),\underline z(\eps)),\qquad (\underline{\underline{x}}(\eps),\underline{\underline{z}}(\eps))^T=P_\eps^+(\underline x(\eps),\underline z(\eps))
\end{equation}
are strictly inside the unit circle. Formula (\ref{P}) yields
\begin{equation}\label{xz__}
   \left(\begin{array}{c}\underline{\underline{x}}(\eps)\\\underline{\underline{z}}(\eps)\end{array}\right)=\left(\begin{array}{c} \underline{x}(\eps)\\ \underline{z}(\eps)-\eps^{1/2}\beta^+h^+(y_0,z_0,0)\sqrt{-M/K}\end{array}\right)\, +\mbox{remaining terms},
\end{equation}
and since
\begin{eqnarray*}
(P_\eps^i)'(x,z)&=& I+\left(\begin{array}{c} 0\\ h^i(y_0,z_0,0)\end{array}\right)\left(\alpha^i,\beta^i+2\eta^i(z-z_0)\right)+\\
&&+\left(\begin{array}{cc}
A^i & 2B^i(z-z_0) \\
0 & 0
\end{array}
\right)\left((\alpha^i,\beta^i)\left(\begin{array}{c} x-x_0\\ z-z_0\end{array}\right)+\eta^i(z-z_0)^2+\gamma^i\eps\right)+\\
&&+\left(\begin{array}{c}A^i(x-x_0)+B^i(z-z_0)^2+C^i\eps \\ 0\end{array}\right)\left(\alpha^i,\beta^i+2\eta^i(z-z_0)\right)+\mbox{remaining terms}, 
\end{eqnarray*}
the composition (\ref{compo}) takes the form
\begin{eqnarray*}
   (P_\eps)'(\underline x(\eps),\underline z(\eps))&=&(\Psi^--\eps^{1/2}(1+\beta^+h^+(y_0,z_0,0))\Phi^-)\left(\Psi^+-\eps^{1/2}\Phi^+\right)+O(\eps^{3/2})=\\
&=&(\Psi^-+\eps^{1/2}\Phi^-)\left(\Psi^+-\eps^{1/2}\Phi^+\right)+O(\eps^{3/2}),
\end{eqnarray*}
where
\begin{eqnarray*}
   \Psi^i&=&I+\left(\begin{array}{cc} 0 & 0 \\
\alpha^i & \beta^i \end{array}\right)h^i(y_0,z_0,0)=\left(\begin{array}{cc} 1 & 0 \\
\alpha^ih^i(y_0,z_0,0)& -1\end{array}\right),\\ 
\Phi^i&=&\left(\begin{array}{cc} A^i\beta^i & 0 \\ 0 & 2h^i(y_0,z_0,0)\eta^i\end{array}\right)\sqrt{-M/K},
\end{eqnarray*}
or, equivalently,
\begin{eqnarray*}
   (P_\eps)'(\underline x(\eps),\underline z(\eps))&=&\Psi^-\Psi^++\eps^{1/2}  (\Phi^-\Psi^+-\Psi^-\Phi^+)+O(\eps)=\\
  &=&\left(\begin{array}{cc} 1 & 0 \\ * & 1\end{array}\right)+\eps^{1/2}\left(\begin{array}{cc} A^-\beta^--A^+\beta^+ & 0 \\ * & -2h^-(y_0,z_0,0)\eta^-+2h^+(y_0,z_0,0)\eta^+\end{array}\right),
\end{eqnarray*}
where *-symbols stay for some (different) constant, which values don't influence the conclusions. Therefore, the fixed point $(\underline x(\eps),\underline z(\eps))$ is stable, if
$$
   h^+(y_0,z_0,0)A^+-h^-(y_0,z_0,0)A^-<0,\qquad  h^+(y_0,z_0,0)\eta^+-h^-(y_0,z_0,0)\eta^-<0.
$$
Analogously, the fixed point $(\overline x(\eps),\overline z(\eps))$ is stable, if
$$
   h^+(y_0,z_0,0)A^+-h^-(y_0,z_0,0)A^->0,\qquad  h^+(y_0,z_0,0)\eta^+-h^-(y_0,z_0,0)\eta^->0.
$$
Assumption (\ref{stab}) ensures that one of these two sets of inequalities holds.

\vskip0.2cm

\noindent {\bf Step 4:} {\it Verifying that $P^-_\eps$ and $P^+_\eps$ map the points of $L$  from the past to the future.}

\vskip0.2cm

\noindent {\it Case 1:} The points from the neighborhood of 
 $(\underline x(\eps),\underline z(\eps))$ and $(\underline {\underline x}(\eps),\underline {\underline z}(\eps))$, i.e. the case where $h^+(y_0,z_0,0)>0.$ We must check that
 \begin{equation}\label{pp1}  
 T^+(\underline x(\eps),\underline z(\eps)>0\quad {\rm and}\quad  
 T^-(\underline {\underline x}(\eps),\underline {\underline z}(\eps))>0.
 \end{equation}
  Using (\ref{Ti}), (\ref{yields}), (\ref{xz}), and (\ref{xz__}), we have
\begin{eqnarray*}
T^+(\underline x(\eps),\underline z(\eps))&=&  \eps^{1/2}\dfrac{2}{h^+(y_0,z_0,0)}\sqrt{-M/K}
+O(\eps),\\
 T^-(\underline {\underline x}(\eps),\underline {\underline z}(\eps))&=&-\eps^{1/2}\dfrac{2}{h^-(y_0,z_0,0)}\sqrt{-M/K}+O(\eps).
\end{eqnarray*}
Therefore the positivity properties (\ref{pp1}) follow from (\ref{timemap}).

\vskip0.2cm

\noindent {\it Case 2:}  The points from the neighborhood of
$(\overline x(\eps),\overline z(\eps))$ and $(\overline {\overline x}(\eps),\overline {\overline z}(\eps)),$ i.e. the case where $h^-(y_0,z_0,0)>0.$  By analogy with Case~1, one can use formulas  (\ref{Ti}), (\ref{yields}), (\ref{xz}), and an analogue of (\ref{xz__}) in order to verify that 
\begin{eqnarray*}
T^+(\overline x(\eps),\overline z(\eps))&=&  -\eps^{1/2}\dfrac{2}{h^+(y_0,z_0,0)}\sqrt{-M/K}
+O(\eps)>0,\\
 T^-(\overline {\overline x}(\eps),\overline {\overline z}(\eps))&=&\eps^{1/2}\dfrac{2}{h^-(y_0,z_0,0)}\sqrt{-M/K}+O(\eps)>0.
\end{eqnarray*}
under  assumption (\ref{timemap}).

\vskip0.2cm

\noindent {\bf Step 5:} {\it Verifying that $P^-_\eps$ and $P^+_\eps$ act in the subspaces $\{(x,y,z):H(y,z)\le0\}$ and 
$\{(x,y,z):H(y,z)\ge0\}$ respectively.}

\vskip0.2cm

\noindent {\it Case 1:} The trajectory with the initial condition at  
 $(\underline x(\eps),\zeta (\underline z(\eps)),\underline z(\eps))$, i.e. the case where $h^+(y_0,z_0,0)>0.$
 
 \vskip0.2cm
 
 \noindent  a) The map $P_\eps^+$. Note, a vector $v\in\mathbb{R}^3$ with the origin $(x,y,z)\in L$ points towards  $\{(x,y,z):H(y,z)\ge0\}$, if $(0,\nabla H(y,z))v>0$. 
 Therefore, the vector field of the ``+''-subsystem of (\ref{sys1}) points to $\{(x,y,z):H(y,z)>0\}$ at $(\underline x(\eps), \underline \zeta(z(\eps)),\underline z(\eps)),$ if
\begin{equation}\label{tmp1}
   \nabla H(\zeta(\underline z(\eps)),\underline z(\eps))\left(\begin{array}{c}
      g(\underline x(\eps),\zeta(\underline z(\eps)))+ g^+(\eps)\\
      h^+(\zeta(\underline z(\eps)),\underline z(\eps),{\color{black}\eps})\end{array}\right)>0,
\end{equation}
which follows from (\ref{directions1}).

\vskip0.2cm 

\noindent b) The map $P_\eps^-.$ Here we have to check that the vector field of the ``-''-subsystem of (\ref{sys1}) points towards $\{(x,y,z):H(y,z)<0\}$ at $(\underline{\underline{x}}(\eps), \zeta(\underline{\underline{z}}(\eps)),\underline{\underline{z}}(\eps)).$ Equivalently, we have to establish that 
\begin{equation}\label{tmp2}   \nabla H(\zeta(\underline{\underline{z}}(\eps)),\underline{\underline{z}}(\eps))\left(\begin{array}{c}
      g(\underline{\underline{x}}(\eps),\zeta(\underline{\underline{z}}(\eps)))+ g^-(\eps)\\
      h^-(\zeta(\underline{\underline{z}}(\eps)),\underline{\underline{z}}(\eps),{\color{black}\eps})\end{array}\right)<0,
\end{equation}
which follows from (\ref{directions2}).

\vskip0.2cm

\noindent {\it Case 2:} The trajectory with the initial condition at  
 $(\overline x(\eps),\zeta(\overline z(\eps)),\overline z(\eps))$,  i.e. the case where $h^+(y_0,z_0,0)<0.$ Considering $(\overline x(\eps),\overline z(\eps))$ in place of $(\underline x(\eps),\underline z(\eps))$ will just flip the sign in the respective expressions (\ref{tmp1})-(\ref{tmp2}), which validity will still follow from (\ref{directions1})-(\ref{directions2}) because the sign of $h^+(y_0,z_0,0)$ flips as well. Therefore,  
 the vector field of the ``+''-subsystem of (\ref{sys1}) points to $\{(x,y,z):H(y,z)>0\}$ at $(\overline x(\eps), \zeta(\overline z(\eps)),\overline z(\eps))$ and 
 the vector field of the ``-''-subsystem of (\ref{sys1}) points to $\{(x,y,z):H(y,z)<0\}$ at $(\overline{\overline{x}}(\eps), \zeta(\overline{\overline{z}}(\eps)),\overline{\overline{z}}(\eps)),$ if conditions (\ref{directions1})-(\ref{directions2}) hold.

\vskip0.2cm

\noindent {\bf Auxiliary relations:} 

\begin{eqnarray*}
\bar k^i&=&\bar{\bar k}^i  \cdot(g(x_0,y_0)+ g^i(0)),\\
\eta^i&=&\dfrac{\bar\eta^i}{h^i(y_0,z_0,0)},\\
\alpha^i&=&\dfrac{\bar\alpha^i}{h^i(y_0,z_0,0)},\\
\gamma^i&=&\dfrac{\bar\gamma^i}{h^i(y_0,z_0,0)},\\
     C^i&=& \dfrac{1}{2}f_y(x_0,y_0)(g(x_0,y_0)+g^i(0))\gamma^i.
\end{eqnarray*}

\begin{eqnarray*}
     B^i&=&\dfrac{1}{2}f'_y(x_0,y_0)(g(x_0,y_0)+g^i(0))\eta^i+\left(\dfrac{1}{6}f'_x(x_0,y_0)g'_y(x_0,y_0)(g(x_0,y_0)+g^i(0))+\right.\\
     && \left.+\dfrac{1}{6}f{}'_y{}'_y(x_0,y_0)(g(x_0,y_0)+g^i(0))^2+\dfrac{1}{6}f'_y(x_0,y_0)g'_y(x_0,y_0)(g(x_0,y_0)+g^i(0))\right)(\beta^i)^2+\\
     &&+\left(\dfrac{1}{2}f_x(x_0,y_0)f_y(x_0,y_0)+\dfrac{1}{2}f_{yy}(x_0,y_0)(g(x_0,y_0)+g^i(0))+\dfrac{1}{2}f_y(x_0,y_0)g_y(x_0,y_0)\right)\dfrac{b}{a+b}\beta^i+\\
&&+\dfrac{1}{2}f{}'_y{}'_y(x_0,y_0)\left(\dfrac{b}{a+b}\right)^2,
\end{eqnarray*}

\noindent The proof of the theorem is complete.

\vskip0.5cm

\section{Numerical computations}\label{sec-comp} In the application, the critical temperature parameter  $T_c^{-}$ is an important bifurcation parameter, below, we show results of the simulations as the system is perturbed by varying this parameter. The following parameters are used:
$$
  b=1.75,\ a=1.05,\ b_0=1.5,\ b_1=5,\ T^-_c=-10,\ \bar T^-=1,\ \bar T^+=0,\ \bar b_0=0,\ \bar b_1=0.
$$
Using Mathematica software we conclude that conditions (\ref{4})-(\ref{as1}) hold for system (\ref{sys0}), if either
\begin{align}
&w = 521.195,\ \eta = -4.50603,\  \xi = -7.80965,\ \ T_c^+=-12.168, \ \ {\rm or}\\ 
&w = -17.2964,\ \eta = 0.244733,\ \xi = -0.208427,\ \ T_c^+=-10.2974, \ \ {\rm or} \\  
&w = 5.08105,\ \eta = 0.948796,\ \xi = 0.918074, \ \ T_c^+=-10.0202. \label{sim3}
\end{align}

\begin{figure}[h]\center
\noindent\includegraphics[scale=0.8]{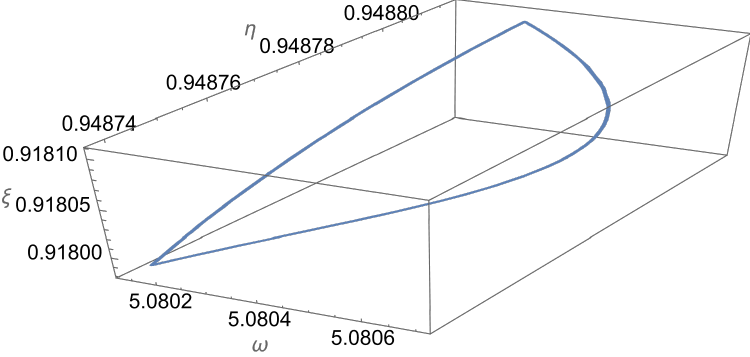}
\caption{Simulation of the attracting cycle of system (\ref{sys0}) with the parameters of Section~\ref{sec-comp} and $\eps=0.001.$} \label{fig-eps001}
\end{figure}

Considering the physically relevant case \eqref{sim3}, with $T^+_c=-10.0202$ and  $\eps=10^{-3},$ formula (\ref{Tformula}) yields a period of $T=0.149956$, while simulations show  $T\approx0.1528.$  Taking $\eps=10^{-4},$ formula (\ref{Tformula}) returns  $T=0.047202$ and simulations show $T\approx0.04795.$ Simulation for $\eps=0.001$ is shown in Fig.~\ref{fig-eps001} 

\section{Closing Remarks}
We have shown that the non-smooth glacial flip flop system \eqref{eqn_wetaxidot} admits a fold fold singularity, and as the critical temperature parameter in the model is varied, the system shows a bifurcation of a limit cycle. Simulations in \cite{widiasih2024mid} show that, as the time dependent parameters \eqref{eqn_Q} and \eqref{eqn_s2} of the orbital forcing are incorporated, varying the same critical temperature parameter produces the Mid Pleistocene transition in the model. A rigorous analysis of such non-autonomous system would be of interest. 


\bibliographystyle{plain}
\bibliography{MWreferences}

\end{document}